\documentclass{article}
\usepackage{amsmath}
\usepackage{amssymb}
\usepackage{latexsym}
\usepackage{amsthm}
\usepackage{mathrsfs}
\usepackage{amsfonts}
\usepackage{amssymb}
\usepackage{graphicx}
\usepackage{latexsym,bm}

\newtheorem{thm}{Theorem}[section]
\newtheorem{lemma}[thm]{Lemma}

\newtheorem{conj}[thm]{Conjecture}
\newtheorem{defi}[thm]{Definition}
\newtheorem{definition}[thm]{Definition}

\newcommand{\JCTB}{{\it J. Combin. Theory Ser. B} }

\begin{document}

\title{Cycle Double Covers and Semi-Kotzig Frame}
\author{Dong Ye and Cun-Quan Zhang\\
{\small Department of Mathematics, West Virginia University,
Morgantown, WV 26506-6310}\\
{\small Emails: dye@math.wvu.edu; cqzhang@math.wvu.edu}}

\maketitle

\begin{abstract}
Let $H$ be a cubic graph admitting a $3$-edge-coloring $c: E(H)\to
\mathbb Z_3$ such that the edges colored by $0$ and $\mu\in\{1,2\}$
induce a Hamilton circuit of $H$ and the edges colored by $1$ and
$2$ induce a 2-factor $F$. The graph $H$ is semi-Kotzig if switching
colors of edges in any even subgraph
 of $F$
 yields a new 3-edge-coloring of $H$ having the
same property as $c$. A spanning subgraph $H$ of a cubic graph $G$
is called a {\em semi-Kotzig frame} if the contracted graph $G/H$ is
even and every non-circuit component of $H$ is a subdivision of a
semi-Kotzig graph.

In this paper, we show that a cubic graph $G$ has a circuit double
cover if it has a semi-Kotzig frame with at most one non-circuit
component. Our result generalizes some results of Goddyn (1988), and
H\"{a}ggkvist and Markstr\"{o}m  [J. Combin. Theory Ser. B
(2006)].
\end{abstract}


\section{Introduction}

Let $G$ be a graph with vertex set $V(G)$ and edge set $E(G)$. A
{\em circuit} of $G$ is a connected 2-regular subgraph. A subgraph
of $G$ is {\em even} if every vertex
 is of even degree. An even subgraph of $G$ is also called a {\em
cycle} in the literatures dealing with cycle covers of graphs
\cite{Jaeger1985} \cite{Jackson1993} \cite{CQ}. Every even graph has
a circuit decomposition. A set $\mathscr C$ of even-subgraphs of $G$
is {\em an even-subgraph double cover} (cycle double cover) if each
edge of $G$ is contained by precisely two even-subgraphs in
$\mathscr C$. The Circuit Double Cover Conjecture was made
independently by Szekeres \cite{GS} and Seymour \cite{PDS}.

\begin{conj}[Szekeres \cite{GS} and Seymour \cite{PDS}]\label{CDCC}
Every bridgeless graph $G$ has a circuit double cover.
\end{conj}

It suffices to show that the Circuit Double Cover Conjecture holds
for bridgeless cubic graphs 
 \cite{Jaeger1985}. The Circuit Double Cover
Conjecture has been verified for several classes of graphs; for
example, cubic graphs with Hamilton paths \cite{Tarsi1986} (also see
\cite{LG2}), cubic graphs with oddness two \cite{HK} and four
\cite{H,HM}, Petersen-minor-free graphs \cite{AGZ}.

A cubic graph $H$ is a {\em spanning minor} of a cubic graph $G$ if
some subdivision of $H$ is a spanning subgraph of $G$. In
\cite{LG1}, Goddyn showed that a cubic graph $G$ has a circuit
double cover if it contains the Petersen graph as a spanning minor.
Goddyn's result is further improved by H\"{a}ggkvist and
Markstr\"{o}m \cite{HMII} who showed that a cubic graph $G$ has a
circuit double cover if it contains a $2$-connected
simple cubic graph with no more than 10 vertices as a
spanning minor.

A {\em Kotzig graph} is a cubic graph $H$ with a 3-edge-coloring $c:
E(G)\to \mathbb Z_3$ such that $c^{-1}(\alpha)\cup c^{-1}(\beta)$
induces a Hamilton circuit of $G$ for every pair $\alpha,\beta\in
\mathbb Z_3$. The family of all Kotzig graphs is denoted by ${\cal
K}$.

\begin{figure}[!hbtp]\refstepcounter{figure}
\label{FIG: Kotzig}
\begin{center}
\includegraphics{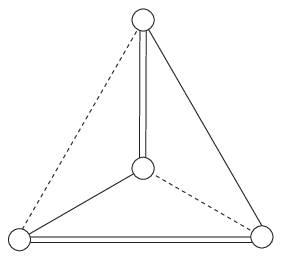}\\
{Figure \ref{FIG: Kotzig}: The Kotzig graph $K_4$. }
\end{center}
\end{figure}

\begin{thm}[Goddyn \cite{LG1}, H\"{a}ggkvist and  Markstr\"{o}m \cite{HMI}]
\label{TH: Kotzig} If a cubic graph $G$ contains a Kotzig graph as a
spanning minor, then $G$ has a $6$-even-subgraph double cover.
\end{thm}

By Theorem \ref{TH: Kotzig}, any cubic graph $G$ containing some
member of ${\cal K}$ as a spanning minor has a circuit double cover.
However, we do not know yet whether every
$3$-connected
  cubic graph contains a member of ${\cal K}$ as a
spanning minor (Conjecture \ref{Kot-SM:Conj1-3}).

According to their observations \cite{HMI,HMII}, H\"{a}ggkvist and
Markstr\"{o}m proposed the following conjectures.

\begin{conj}[H\"{a}ggkvist and Markstr\"{o}m,
\cite{HMI}]
\label{Kot-SM:Conj1-3} Every $3$-connected cubic graph
contains a Kotzig graph as a spanning minor\footnote{
It is pointed out in 
\cite{HoffmannOstenhof2011} that the $3$-edge-connectivity is not enough for the existence of such spanning minor, and he suggested that an extra requirement of cyclical $4$-edge-connectivity is necessary.}.
\end{conj}

H\"{a}ggkvist and Markstr\"{o}m \cite{HMI} proposed another
conjecture (Conjecture~\ref{3ECC}) in a more general form. We will
discuss this conjecture in the last section (Remark).

\begin{figure}[!hbtp]\refstepcounter{figure}
\label{FIG: Kotzig-inclusion}
\begin{center}
\includegraphics{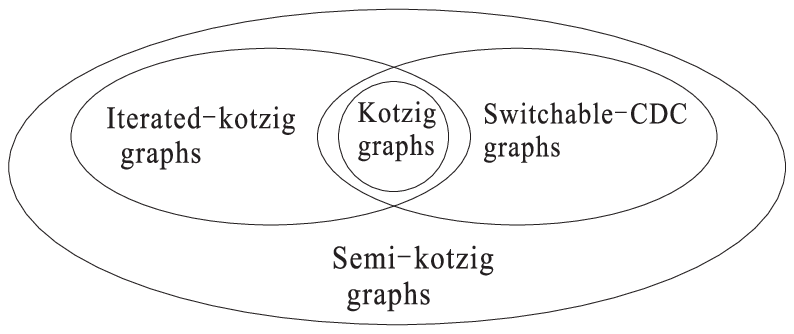}\\
{Figure \ref{FIG: Kotzig-inclusion}: The inclusion relations for
these four families: Kotzig graphs, iterated-Kotzig graphs,
switchable CDC graphs, Semi-Kotzig graphs. }
\end{center}
\end{figure}

One of approaches to the CDC conjecture is to find a sup-family
${\cal X}$  of ${\cal K}$ such that every bridgeless cubic graph
containing a member of ${\cal X}$ as a spanning minor has a CDC.
Following this direction of approach, Goddyn \cite{LG1},
H\"{a}ggkvist and Markstr\"{o}m \cite{HMI} introduce some
sup-families of ${\cal K}$, named iterated-Kotzig graphs,
switchable-CDC graphs and semi-Kotzig graphs. They will be defined
in next subsections and their relations are shown in
Figure~\ref{FIG: Kotzig-inclusion}.

\subsection*{Iterated-Kotzig graphs}

\begin{defi}
\label{DEF: Iterated}
 {\em An {\em iterated-Kotzig graph} $H$ is a cubic graph
constructed as following \cite{HMI}: Let $\mathcal K_0$ be a set of
Kotzig graphs with a 3-edge-coloring $c: E(G)\to \mathbb Z_3$; A
cubic graph $H\in \mathcal K_{i+1}$ can be constructed from a graph
$H_i\in \mathcal K_i$ and a graph $H_0\in \mathcal K_0$ by deleting
one edge colored by 0 from each of them and joining the two vertices
of degree two in $H_0$ to the two vertices of degree two in $H_i$,
respectively (the two new edges will be colored by 0).}
\end{defi}

\begin{figure}[!hbtp]\refstepcounter{figure}
\label{FIG: iterated-Kotzig}
\begin{center}
\includegraphics{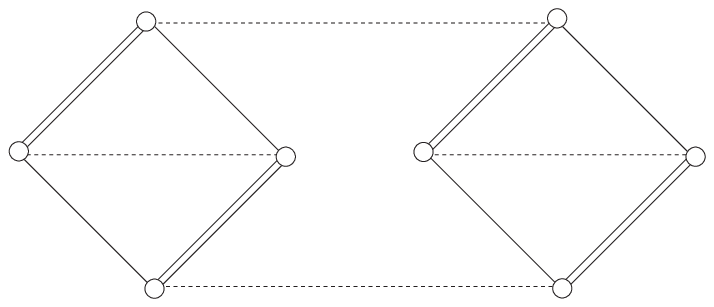}\\
{Figure \ref{FIG: iterated-Kotzig}: An iterated-Kotzig graph
generated from two $K_4$'s. }
\end{center}
\end{figure}

\begin{thm}[H\"{a}ggkvist and Markstr\"{o}m, \cite{HMI}]
\label{TH: iterated-one} If a cubic graph $G$ contains an iterated
Kotzig graph as a spanning minor, then $G$ has a $6$-even-subgraph
double cover.
\end{thm}


\subsection*{Semi-Kotzig graphs and switchable-CDC graphs}
\begin{defi}
\label{DEF: semi-K} {\upshape Let $G$ be a cubic graph with a
3-edge-coloring $c: E(G)\to \mathbb Z_3$ and the following property
\begin{center}
($*$) \hspace{0.3cm} edges in colors 0 and $\mu$ ($\mu\in \{1,2\}$)
induce a Hamilton circuit.
\end{center}
Let $F$ be the even 2-factor induced by edges in colors 1 and 2. If,
for every even subgraph $S\subseteq F$, switching colors 1 and 2 of
the edges of $S$
 yields a new 3-edge-coloring having the property $(*)$, then
each of these $2^{t-1}$ 3-edge-coloring is called a {\em semi-Kotzig
coloring} where $t$ is the number of components of $F$. A cubic
graph $G$ with a semi-Kotzig coloring is called a {\em semi-Kotzig
graph}. If $F$ has at most two components ($t \leq 2$), then $G$ is
said to be a {\em switchable-CDC graph} (defined in \cite{HMI}).}
\end{defi}

\begin{figure}[!hbtp]\refstepcounter{figure}
\label{FIG: Semi-Kotzig}
\begin{center}
\includegraphics{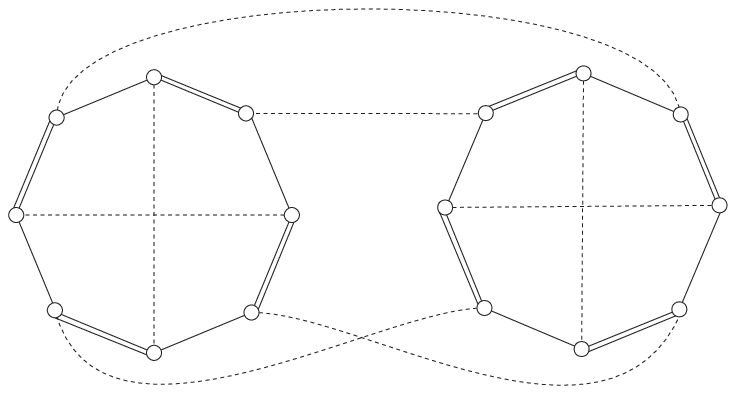}\\
{Figure \ref{FIG: Semi-Kotzig}: A semi-Kotzig graph. }
\end{center}
\end{figure}

\begin{thm}[H\"{a}ggkvist and Markstr\"{o}m,
\cite{HMI}] \label{TH: 2-semi-Kotzig} If a cubic graph $G$ contains
a switchable-CDC graph as a spanning minor, then $G$ has
  a $6$-even-subgraph double cover.
\end{thm}

An iterated-Kotzig graph has a semi-Kotzig coloring and hence is a
semi-Kotzig graph. But a semi-Kotzig graph is not necessary an
iterated-Kotzig graph. For example, the semi-Kotzig graph in Figure
\ref{FIG: Semi-Kotzig} is not an iterated-Kotzig graph. Hence we
have the following relations (also see Figure~\ref{FIG:
Kotzig-inclusion})
\begin{equation}
\mbox{Kotzig} ~ \subset ~ \mbox{ Iterated-Kotzig} ~ \subset ~ \mbox{
Semi-Kotzig}; \label{EQ: Kotzig inclusion}
\end{equation}
\begin{equation}
\mbox{Kotzig} ~ \subset ~ \mbox{ Switchable-CDC} ~ \subset ~ \mbox{
Semi-Kotzig}.\label{EQ: Switchable-CDC inclusion}
\end{equation}

The following theorem was announced in \cite{LG1} with an outline of
proof.

\begin{thm}[Goddyn, \cite{LG1}]\label{TH: semi-Kotzig}
If a cubic graph $G$ contains a semi-Kotzig graph as a spanning
minor, then $G$ has a $6$-even-subgraph double cover.
\end{thm}

The main theorem (Theorem~\ref{TH: SemiKotzig-disconnected})
 of the paper strengthens all those early results
(Theorems \ref{TH: Kotzig}, \ref{TH: iterated-one}, \ref{TH:
2-semi-Kotzig} and \ref{TH: semi-Kotzig}).

\subsection*{Kotzig frame, semi-Kotzig frame}

A $2$-factor $F$ of a cubic graph is {\em even} if every component
of $F$ is of even length. If a cubic graph $G$ has an even
$2$-factor, then the graph $G$ has many nice properties: {\em $G$ is
$3$-edge-colorable, $G$ has a circuit double cover and strong
circuit double cover, etc.}

 The
following concepts were introduced in \cite{HMI} as a generalization
of even $2$-factors.

\begin{defi} {\upshape Let $G$ be a cubic graph.  A spanning subgraph $H$ of $G$ is called
a {\em frame} of $G$ if the contracted graph $G/H$ is an even
graph.}
\end{defi}

For a subgraph $H$ of $G$, the {\em suppressed graph} $\overline{H}$
of $H$ is the graph obtained from $H$ by suppressing all degree $2$
vertices.

\begin{defi}
{\em Let $G$ be a cubic graph.  A frame $H$ of $G$ is called a {\em
Kotzig frame} (or {\em iterated-Kotzig frame}, or {\em
switchable-CDC frame}, or {\em semi-Kotzig frame}) of $G$ if, for
each non-circuit component $H_i$ of $H$, the suppressed graph
$\overline{H_i}$ is a Kotzig graph (or an iterated-Kotzig graph, or
a switchable-CDC graph, or a semi-Kotzig graph, respectively).}
\end{defi}

Similar to the relations described in (\ref{EQ: Kotzig inclusion})
and (\ref{EQ: Switchable-CDC inclusion}), we have the same relations
between those frames:
\[\mbox{Kotzig frame} ~ \subset ~ \mbox{ Iterated-Kotzig\ frame} ~ \subset ~ \mbox{ semi-Kotzig\ frame};\]
\[\mbox{Kotzig frame} ~ \subset ~ \mbox{ Switchable-CDC\ frame} ~ \subset ~ \mbox{ semi-Kotzig\ frame}.\]

\begin{thm}[H\"aggkvist and Markstr\"om, \cite{HMI}]\label{TH: Kotzig-disconnected}
Let $G$ be a bridgeless cubic  graph $G$. If $G$ contains a Kotzig
frame with at most one non-circuit component,
  then
$G$ has a $6$-even-subgraph double cover.
\end{thm}

According to their observations, they further make the following
conjecture.

\begin{conj}[H\"aggkvist and Markstr\"om, \cite{HMI}]\label{CONJ: Kotzig-frame}
Every bridgeless cubic graph with a Kotzig frame
  has a $6$-even-subgraph double cover.
\end{conj}

The following theorem provides a partial solution to Conjecture
\ref{CONJ: Kotzig-frame}.

\begin{thm}[Zhang and Zhang, \cite{ZZ}]
\label{TH: Kotzig - tree} Let $G$ be a bridgeless cubic graph. If
$G$ contains a Kotzig frame $H$ such that $G/H$ is a tree if
parallel edges are identified as a single edge,
  then $G$ has a $6$-even-subgraph double cover.
\end{thm}

We conjecture that the result in Conjecture \ref{CONJ: Kotzig-frame}
still holds if a
Kotzig frame is replaced by a semi-Kotzig frame.

\begin{conj}\label{CONJ: Semi-Kotzig-frame}
Every bridgeless cubic graph with a semi-Kotzig frame
 has a $6$-even-subgraph double cover.
\end{conj}

H\"aggkvist and Markstr\"om showed Conjecture \ref{CONJ:
Semi-Kotzig-frame} holds for iterated-kotzig frames and
switchable-CDC frames with at most one non-circuit component.

\begin{thm} [H\"aggkvist and Markstr\"om, \cite{HMI}]\label{TH:
iterated-disconnected} Let $G$ be a bridgeless cubic  graph $G$. If
$G$ contains an iterated-Kotzig frame with at most one non-circuit
component, then $G$ has a $6$-even-subgraph double cover.
\end{thm}

\begin{thm} [H\"aggkvist and Markstr\"om, \cite{HMI}]\label{TH:
Switchable-CDC frame} Let $G$ be a bridgeless cubic  graph $G$. If
$G$ contains a switchable-CDC frame with at most one non-circuit
component, then $G$ has a $6$-even-subgraph double cover.
\end{thm}

The following theorem is the main result of the paper, which
verifies that Conjecture \ref{CONJ: Semi-Kotzig-frame} holds if a
semi-Kotzig frame has at most one non-circuit component. Since
Kotzig graphs, iterated-Kotzig graphs are semi-Kotzig graphs but not
vice verse, Theorems \ref{TH: Kotzig}, \ref{TH: iterated-one},
\ref{TH: 2-semi-Kotzig}, \ref{TH: semi-Kotzig}, \ref{TH:
Kotzig-disconnected}, \ref{TH: iterated-disconnected} and \ref{TH:
Switchable-CDC frame} are corollaries of our result. The proof of
the theorem will be given in Section 2.

\begin{thm}
\label{TH: SemiKotzig-disconnected} Let $G$ be a bridgeless cubic
graph. If $G$ contains a semi-Kotzig frame $H$ with at most one
non-circuit component,
  then $G$ has a $6$-even-subgraph double cover.
\end{thm}

\section{Proof of Theorem \ref{TH: SemiKotzig-disconnected}}

The following well-known fact will be applied in the proof of the
main theorem (Theorem~\ref{TH: SemiKotzig-disconnected}).

\begin{lemma}
\label{LE: 3-edge-color, 2-cover} If a cubic graph has an even
$2$-factor $F$, then
  $G$ has a  $3$-even-subgraph
double cover ${\mathscr C}$ such that $F \in {\mathscr C}$.
\end{lemma}

\begin{definition}
{\em Let $H$ be a bridgeless subgraph of a cubic graph $G$.  A
mapping $c: E(H) \to \mathbb Z_3$ is called a {\em parity
$3$-edge-coloring }of $H$ if, for each vertex $v \in H$ and each
$\mu \in \mathbb Z_3$,
$$|c^{-1}(\mu) \cap E(v)| \equiv |E(v) \cap E(H)| ~~ \pmod{2}.$$
It is obvious that if $H$ itself is cubic, then a parity
$3$-edge-coloring is a proper $3$-edge-coloring (traditional
definition).}
\end{definition}

\medskip \noindent
{\bf Preparation of the proof}. Let $H_0$ be the component of $H$
such that $H_0$ is a subdivision of a semi-Kotzig graph and each
$H_i$, $1\le i\le t$, be a circuit component of $H$ of even length.
Let $M=E(G)-E(H)$, and $H^* = H- H_0$.

Given an initial
 semi-Kotzig
coloring $c_0:E(\overline H_0)\to \mathbb Z_3$ of $\overline H_0$,
then $F_0= c_0^{-1}(1) \cup c_0^{-1}(2)$ is a 2-factor of $\overline
H_0$ and $c_0^{-1}(0)\cup c_0^{-1}(\mu)$ is a Hamilton circuit of
$\overline H_0$ for each $\mu\in \{1,2\}$.

The semi-Kotzig coloring $c_0$ of $\overline H_0$ can be considered
as an edge-coloring of $H_0$: each induced path is colored with the
same color as its corresponding edge in $\overline H_0$
  (note, this
edge-coloring of $H_0$ is a parity $3$-edge-coloring, which may not
be a proper $3$-edge-coloring).

\label{P: strategy} The strategy of the proof is to show that $G$
can be covered by three subgraphs $G(0,1), G(0,2)$ and $G(1,2)$ such
that each $G(\alpha, \beta)$ has a 2-even-subgraph cover which
covers the edges of $M \cap E(G(\alpha, \beta))$ twice and the edges
of $E(H) \cap E(G(\alpha, \beta))$ once. In order to prove this, we
are going to show that the three subgraphs $G(\alpha, \beta)$ have
the following properties:

(i) the suppressed cubic graph $\overline{G(\alpha, \beta)}$ is
$3$-edge-colorable (so that Lemma~\ref{LE: 3-edge-color, 2-cover}
can be applied to each of them);

(ii) $c_0^{-1}(\alpha) \cup c_0^{-1}(\beta) \subseteq G(\alpha,
\beta)$ for each pair $\alpha, \beta \in \mathbb Z_3$;

(iii) The even subgraph $H^*$ has a decomposition, $H^*_1$ and
$H^*_2$, each of which is an even subgraph, (here, for technical
reason, let $H^*_0=\emptyset$), such that $H^*_{\alpha} \cup
H^*_{\beta} \subseteq G(\alpha, \beta)$, for each $\{ \alpha, \beta
\} \subset \mathbb Z_3$;

(iv) each $e \in M = E(G) - E(H)$ is contained in precisely one
member of $\{ G(0,1), G(0,2), G(1,2)\}$;

(v) and most important, the subgraph $c^{-1}(\alpha) \cup
c^{-1}(\beta) \cup H^{*}_{\alpha} \cup H^{*}_{\beta}$ in
$G(\alpha,\beta)$ corresponds to an even $2$-factor of $
\overline{G(\alpha, \beta)}$.

Can we decompose $H^*$ and find a partition of $ M = E(G) - E(H)$ to
satisfy (v)?
One may also notice that the
initial semi-Kotzig coloring $c$ may not be appropriate. However,
the color-switchability of the semi-Kotzig component $H_0$ may help
us to achieve the goal. The properties described above in the
strategy will be proved in the following claim.

\medskip
We claim that $G$ has the following property: \medskip
\begin{description}

\item{$(*)$}~ {\sl There is a semi-Kotzig coloring $c_0$ of
$\overline H_0$, a decomposition $\{ H^*_{1}, H^*_2 \}$ of $H^*$ and
a partition $\{ N_{(0,1)}, N_{(0,2)}, N_{(1,2)}\}$ of $M$ such that,
let $C_{(\alpha,\beta)} = c_0^{-1}(\alpha) \cup
c_0^{-1}(\beta)$,

{\upshape(1)}
for each $\mu\in \{1,2\}$, $C_{(0,\mu)}\cup
H^{*}_{\mu}$ corresponds to an even $2$-factor of
$\overline{G(0,\mu)} =\overline{ G[C_{(0,\mu)} \cup H^*_{\mu} \cup
N_{(0,\mu)}]}$, and

{\upshape(2)} $C_{(1,2)}\cup H^*$ corresponds to an even $2$-factor
of $\overline{G(1,2)} =\overline{G[C_{(1,2)} \cup   H^* \cup
N_{(1,2)}]}$.}\end{description}

\medskip
\noindent {\bf Proof of $(*)$.} Let $G$ be a minimum counterexample
to ($*$). Let $c: E(H) \to \mathbb Z_3$ be a parity
$3$-edge-coloring of $H$ such that

(1) the restriction of $c$ on $\overline{H_0}$ is a semi-Kotzig
coloring, and

(2) $E(H^*) \subseteq c^{-1}(1) \cup c^{-1}(2)$ (a set of
mono-colored circuits).\\ Let \[F = c^{-1}(1) \cup c^{-1}(2) = E(H)
- c^{-1}(0).\]

Partition the matching $M$ as follows. For each edge $e=xy \in M$,
$xy \in M_{(\alpha, \beta)}$  ($\alpha \leq \beta$ and $\alpha,
\beta \in \mathbb Z_3$) if $x$ is incident with two $\alpha$-colored
edges and $y$ is incident with two $\beta$-colored edges. So, the
matching $M$ is partitioned into six subsets:
\[M_{(0,0)}, M_{(0,1)}, M_{(0,2)}, M_{(1,1)}, M_{(1,2)}\mbox{ and } M_{(2,2)}.\] Note
that this partition will be adjusted whenever the parity
$3$-edge-coloring is adjusted.

\medskip \noindent
{\bf Claim 1.} {\sl $M_{(0,\mu)}\cap G[V(H_0)]=\emptyset$, for each
$\mu \in \mathbb Z_3$.}

\vspace{0.3cm}

Suppose that $e=xy\in M_{(0,\mu)}$ where $x$ is incident with two
0-colored edges of $H_0$. Then, in the graph $\overline{G-e}$, the
spanning subgraph $H$ retains the same property as itself in $G$.
Since $\overline{G-e}$ is smaller than $G$, $\overline{G-e}$
satisfies ($*$): $\overline H_0$ has a semi-Kotzig coloring $c_0$
and $M-e$ has a partition $\{ N_{(0,1)}, N_{(0,2)}, N_{(1,2)}\}$ and
$H^*$ has a decomposition $\{ H_1^*, H_2^*\}$. In the semi-Kotzig
coloring $c_0$, without loss of generality, assume that $y$
subdivides a 1-colored edge of $\overline H_0$. For the graph $G$,
add $e$ into $N_{(0,1)}$. This revised partition $\{N_{(0,1)},
N_{(0,2)}, N_{(1,2)}\}$ of $M$ and the resulting subgraphs
$G(\alpha,\beta)$ satisfy ($*$). This contradicts that $G$ is a
counterexample.\medskip

Since $c^{-1}(0) \subseteq H_0$ (each component of $H-H_0=H^*$ is
mono-colored by $1$ or $2$), for every edge $e \in M_{(0,\mu)}$
($\mu \in \{ 1,2 \}$), by Claim 1, the edge $e$ has one endvertex
incident with two $0$-colored edges of $H_0$ and its another
endvertex belongs to $V(H-H_0)=V(H^*)$. That is,
$$M_{(0,0)}=\emptyset, ~~\mbox{and}~~
M_{(0,1)} \cup M_{(0,2)} \subseteq E(H_0, H^*).$$

Let \[G'=\overline{G-M_{(0,1)}-M_{(0,2)}}.\] Then $E(G'/F)\subseteq
M_{(1,1)}\cup M_{(1,2)}\cup M_{(2,2)}$.

\medskip

\noindent {\bf Claim 2.} {\sl The graph $G'/F$ is acyclic.}\medskip


Suppose to the contrary that $G'/F$ contains a circuit $Q$
(including loops). In the graph $\overline{G-E(Q)}$, the spanning
subgraph $H$ retains as a semi-Kotzig frame.

Then the smaller graph $\overline{G-E(Q)}$ satisfies $(*)$:
$\overline H_0$ has a semi-Kotzig coloring $c_0$, and  $M-E(Q)$ has
a a partition $\{ N_{(0,1)}, N_{(0,2)}, N_{(1,2)}\}$, and $H^*$ has
a decomposition $\{H_1^*, H_2^*\}$.
  So add all edges of $E(Q)$ into $N_{(1,2)}$.
This revised partition $\{N_{(0,1)}, N_{(0,2)}, N_{(1,2)}\}$ of $M$
and its resulting subgraphs $G(\alpha,\beta)$ also satisfy $(*)$
since $C_{(1,2)}\cup H^*$ corresponds to an even $2$-factor of
$\overline{G(1,2)} =\overline{G[C_{(1,2)} \cup H^* \cup
N_{(1,2)}]}$. This is a contradiction. So Claim 2 follows. \medskip

By Claim 2, each component $T$ of $G'/F$ is a tree. Along the tree
$T$, we can modify the parity 3-edge-coloring $c$ of $H$ as follows:

\begin{description}
\item {$(**)$} properly switch colors for some circuits in $ F$
so that every edge of $T$ is incident with four same colored edges.
\end{description}

Note that Rule ($**$) is feasible by Claim 2 since $G'/F$ is
acyclic. Furthermore, under the modified parity 3-edge-coloring $c$,
$M_{(1,2)}=\emptyset$. So
\[M=M_{(0,1)}\cup M_{(0,2)}\cup M_{(1,1)}\cup M_{(2,2)}.\]
The colors of all $H_i$'s $(i\ge 1)$ give a decomposition $\{H_1^*,
H_2^*\}$ of $H^*$ where $H_{\mu}^*$ consists of all circuits of
$H^*$ mono-colored by $\mu$ for $\mu=1$ and 2.

Let \[G''=G/H\]
where $E(G'')=M$. Then $G''$ is even since $H$ is a
frame. For a vertex $w$ of $G''$ corresponding to a component $H_i$
with $i\ge 1$, there is a $\mu \in \{ 1,2 \}$ such that
  all edges incident with $w$
belong to $ M_{(0,\mu)}\cup M_{(\mu,\mu)}$. Define
$$N_{(0,\mu)}=M_{(0,\mu)}\cup M_{(\mu,\mu)}$$ for each
$\mu\in\{1,2\}$, and
\[N_{(1,2)}=M_{(1,2)}=\emptyset.\]
Hence, a
vertex of $G''$ corresponding to $H_i$ with $i\ge 1$ either has
degree in $G''[N_{(0,\mu)}]$
 the same as its degree in $G''$ or
has degree zero (by Rule $(**)$).
 So every vertex of $G''[N_{(0,\mu)}]$ which is
different from the vertex corresponding to $H_0$ has even degree.
Since every graph has even number of odd-degree vertices, it follows
that $G''[N_{(0,\mu)}]$ is an even subgraph.

For each $\mu \in \{ 1,2 \}$, let $G(0,\mu) = N_{(0,\mu)} \cup
(c^{-1}(0) \cup c^{-1}(\mu))$. Since $G''[N_{(0,\mu)}]$ is an even
subgraph of $G''$, the even subgraph $c^{-1}(0) \cup c^{-1}(\mu)$
corresponds to an even $2$-factor of $G(0,\mu)$. And let $G(1,2) = F
=c^{-1}(1)\cup c^{-1}(2)$ (here,
$N_{(1,2)}=\emptyset$).
  So $G$ has the
property ($*$), a contradiction. This completes the proof of ($*$).
\qed

\medskip \noindent
{\bf Proof of Theorem~\ref{TH: SemiKotzig-disconnected}.} Let $G$ be
a graph with a semi-Kotzig frame. Then $G$ satisfies $(*)$ and
therefore is covered by three subgraphs $G(\alpha, \beta)$
($\alpha,\beta\in \mathbb Z_3$ and $\alpha<\beta$) as stated in
$(*)$.

Applying Lemma \ref{LE: 3-edge-color, 2-cover} to the three graphs
$\overline{G(\alpha, \beta)}$, each $G(0,\mu)$ has a 2-even-subgraph
cover $\mathscr C_{(0,\mu)}$ which covers the edges of
$C_{(0,\mu)}\cup H^*_{\mu}$ once and the edges in $N_{(0,\mu)}$
twice, and $G(1,2)$ has a 2-even-subgraph cover $\mathscr C_{(1,2)}$
which covers the edges of $C_{(1,2)}\cup H^*$ once and the edges in
$N_{(1,2)}$ twice. So $\bigcup \mathscr C_{(\alpha, \beta)}$
  is a $6$-even-subgraph double cover of $G$.
  This completes the proof.\qed

\vspace{.3cm}

\noindent{\bf Remark.} In \cite{HMI}, H\"{a}ggkvist and
Markstr\"{o}m proposed another conjecture which strengthens Theorems
\ref{TH: Kotzig}, \ref{TH: iterated-one} and \ref{TH: semi-Kotzig}
as follows.

\begin{conj}[H\"{a}ggkvist and Markstr\"{o}m,
\cite{HMI}]\label{3ECC} If a cubic bridgeless graph contains a connected
$3$-edge-colorable cubic graph as a spanning minor,
 then $G$ has a $6$-even-subgraph double cover
\end{conj}

In fact, Conjecture~\ref{3ECC} is equivalent to that every
bridgeless cubic graph has a $6$-even-subgraph double cover. It can be
shown that the condition in Conjecture~\ref{3ECC} is true for all
cyclically $4$-edge-connected cubic graphs.

Consider a cyclically $4$-edge-connected cubic graph $G$ since a
smallest counterexample to
 the $6$-even-subgraph double cover problem is
cyclically $4$-edge-connected and cubic.
 By the
Matching Polytop Theorem of Edmonds \cite{JE}, $G$ has a $2$-factor
$F$ such that $G/F$ is 4-edge-connected. By Tutte \& Nash-Williams
Theorem (\cite{NW, Tutte61}), $G/F$ contains two edge-disjoint
spanning trees $T_1$ and $T_2$. By a theorem of Itai and Rodeh
(\cite{Itai1978}), $T_1$ contains a parity subgraph $P$ of $G/F$.
 After suppressing all degree two vertices of $G-P$, the graph
$\overline{G-P}$ is $3$-edge-colorable and connected since $G/F-P$
is even and $T_2\subset G/F-P$. So every cyclically
$4$-edge-connected cubic graph does contain a connected
$3$-edge-colorable cubic graph as a spanning minor.


\end{document}